\documentclass{amsart}
\usepackage{amssymb}
\usepackage{amsmath}
\usepackage{textcomp}
\usepackage{color}
\usepackage{graphicx}
\usepackage{pict2e}

\newtheorem{theorem}{Theorem}[section]
\newtheorem{lemma}[theorem]{Lemma}

\theoremstyle{definition}

\newtheorem{corollary}[theorem]{Corollary}
\theoremstyle{remark}

\numberwithin{equation}{section}

\begin{document}

\title{Norm On The Variable Mixed Space $ \ell^{q(\cdot)}(L^{p(\cdot)}) $}



\date{}
\author{Reza Roohi Seraji}
\address{
Department of Mathematics, \\
Institute for Advanced Studies in Basic Sciences (IASBS),\\ 
Zanjan 45137-66731, Iran
}
\email{rroohi@iasbs.ac.ir}
\dedicatory{}

\commby{}


\maketitle

\begin{abstract}
This paper aims to establish the norm properties of the variable mixed space $ \ell^{q(\cdot)}(L^{p(\cdot)}) $ when $ 1<q_-,p_-,q_+,p_+<\infty $. In this way, we address the open problem raised by Almeida and H\"{a}st\"{o}.
\end{abstract}

\section{Preliminaries}

Denote by $\mathcal{P}_0$ the set of measurable functions $p(\cdot): \mathbb{R}^n \rightarrow [c,\infty]$, where $c>0$, and by $\mathcal{P}$ the subset of $\mathcal{P}_0$ such that the support of its elements is contained in $[1,\infty]$. If $p(\cdot), q(\cdot) \in \mathcal{P}_0$, then the norm of the space $\ell^{q(\cdot)} (L^{p(\cdot)}(\mathbb{R}^n))$ naturally is defined as follows is a quasi-norm,
\begin{align}\label{def:norm}
\|f\|_{{\ell}^{q(\cdot)}(L^{p(\cdot)})}&= \inf\left\{\mu>0 ~\Big|~ \varrho_{{\ell}^{q(\cdot)}(L^{p(\cdot)})}\left(\mu^{-1}f\right)\le 1\right\}
\end{align}
where 
\begin{align}\label{def:modular}
\varrho_{\ell^{q(\cdot)}(L^{p(\cdot)})}(f)=\sum_{\nu=1}^{\infty}\inf\{\lambda_\nu>0~|~\varrho_{p(\cdot)}(\frac{f_\nu}{\lambda_{\nu}^{\frac{1}{q(\cdot)}}})\leq 1\}.
\end{align}

\section{Norm Of Space}
\begin{theorem}\label{theo:SctrictlyConvex}
\cite[Corollary 2.8]{GR} Let $ 1<p_-,q_-,p_+,q_+<\infty $, then $ \ell^{q(\cdot)}(L^{p(\cdot)}) $ is strictly convex.
\end{theorem}

\begin{lemma}\label{lem:StrictInequality}
Let $ \ell^{q(\cdot)}(L^{p(\cdot)}) $ is strictly convex and for $ \nu\in\mathbb{N} $, we have
\begin{align*}
1<\varrho_{p(\cdot)}(\frac{f_\nu}{\zeta_\nu^{\frac{1}{q(\cdot)}}(\mu_1+\mu_2)}+\frac{g_\nu}{\zeta_\nu^{\frac{1}{q(\cdot)}}(\mu_1+\mu_2)})
\end{align*}
then
\begin{align*}
\exists r\in(0,\mu_2-\mu_1]~:~1\lneqq\varrho_{p(\cdot)}(\frac{f_\nu}{\zeta_\nu^{\frac{1}{q(\cdot)}}(\mu_1+\mu_2-r)}+\frac{g_\nu}{\zeta_\nu^{\frac{1}{q(\cdot)}}(\mu_1+\mu_2+r)})
\end{align*}
where $ \mu_1\geq\|f\|_{\ell^{q(\cdot)}(L^{p(\cdot)})},~ \mu_2\geq\|g\|_{\ell^{q(\cdot)}(L^{p(\cdot)})} $ and $ \mu_1\leq\mu_2 $.
\end{lemma}
\begin{proof}
According to strictly convexity, for $ f\slash\|f\|_{\ell^{q(\cdot)}(L^{p(\cdot)})} $ and $ g\slash\|g\|_{\ell^{q(\cdot)}(L^{p(\cdot)})} $ we have
\begin{align}\label{eq:StrictlyConvex}
\varrho_{\ell^{q(\cdot)}(L^{p(\cdot)})}(\frac{f}{2\|f\|_{\ell^{q(\cdot)}(L^{p(\cdot)})}}+\frac{g}{2\|g\|_{\ell^{q(\cdot)}(L^{p(\cdot)})}})<1,
\end{align}
which turns into,  there exists $ \{\zeta_\nu\}_\nu $ such that $ \sum_{\nu=1}^{\infty}\zeta_\nu\leq 1 $ and for every $ \nu $ we have
\begin{align}\label{eq:StriclyConvex2}
\varrho_{p(\cdot)}(\frac{f_\nu}{\zeta_\nu^{\frac{1}{q(\cdot)}} 2\mu_1}+\frac{g_\nu}{\zeta_\nu^{\frac{1}{q(\cdot)}} 2\mu_2})<1.
\end{align}
Now we are going to show
\begin{align*}
\exists r\in(0,\mu_2-\mu_1]~:~1=\varrho_{p(\cdot)}(\frac{f_\nu}{\zeta_\nu^{\frac{1}{q(\cdot)}}(\mu_1+\mu_2+r)}+\frac{g_\nu}{\zeta_\nu^{\frac{1}{q(\cdot)}}(\mu_1+\mu_2+r)}).
\end{align*}
By reductio and absurdum we have
\begin{align*}
\forall r\in(0,\mu_2-\mu_1]~:~1<\varrho_{p(\cdot)}(\frac{f_\nu}{\zeta_\nu^{\frac{1}{q(\cdot)}}(\mu_1+\mu_2+r)}+\frac{g_\nu}{\zeta_\nu^{\frac{1}{q(\cdot)}}(\mu_1+\mu_2+r)}).
\end{align*}
Let $ \mu_1\leq\mu_2,~ r:=\mu_2-\mu_1 $, and taking into account the relation \eqref{eq:StriclyConvex2} we have
\begin{align*}
1&<\varrho_{p(\cdot)}(\frac{f_\nu}{\zeta_\nu^{\frac{1}{q(\cdot)}}2\mu_2}+\frac{f_\nu+g_\nu}{\zeta_\nu^{\frac{1}{q(\cdot)}}2\mu_2})
\\
&\leq\varrho_{p(\cdot)}(\frac{f_\nu}{\zeta_\nu^{\frac{1}{q(\cdot)}}2\mu_1}+\frac{f_\nu+g_\nu}{\zeta_\nu^{\frac{1}{q(\cdot)}}2\mu_2})
<1,
\end{align*}
which is a contradiction. Therefore,
\begin{align*}
\exists r\in(0,\mu_2-\mu_1]~:~1&=\varrho_{p(\cdot)}(\frac{f_\nu}{\zeta_\nu^{\frac{1}{q(\cdot)}}(\mu_1+\mu_2+r)}+\frac{g_\nu}{\zeta_\nu^{\frac{1}{q(\cdot)}}(\mu_1+\mu_2+r)})
\\
&\lneqq\varrho_{p(\cdot)}(\frac{f_\nu}{\zeta_\nu^{\frac{1}{q(\cdot)}}(\mu_1+\mu_2-r)}+\frac{g_\nu}{\zeta_\nu^{\frac{1}{q(\cdot)}}(\mu_1+\mu_2+r)})
\end{align*}
which is the desired result.
\end{proof}
\begin{corollary}\label{coro:norm}
Let $ 1<p_-,q_-,p_+,q_+<\infty $, then the definition in \eqref{def:norm} defines a norm in $ \ell^{q(\cdot)}(L^{p(\cdot)}) $.
\end{corollary}
\begin{proof}
It suffices to investigate the triangle inequality i.e. we need to show that for every $ f,g\in\ell^{q(\cdot)}(L^{p(\cdot)}) $ we have
\begin{align*}
\|f+g\|_{\ell^{q(\cdot)}(L^{p(\cdot)})}\leq \|f\|_{\ell^{q(\cdot)}(L^{p(\cdot)})}+\|g\|_{\ell^{q(\cdot)}(L^{p(\cdot)})}.
\end{align*}
Let $ \mu_1\geq \|f\|_{\ell^{q(\cdot)}(L^{p(\cdot)})} $ and $ \mu_2\geq\|g\|_{\ell^{q(\cdot)}(L^{p(\cdot)})} $, then according to definition \eqref{def:norm}, we need to show
\begin{align*}
\varrho_{\ell^{q(\cdot)}(L^{p(\cdot)})}(\frac{f+g}{\mu_1+\mu_2})\leq 1.
\end{align*}
For $ f\slash\|f\|_{\ell^{q(\cdot)}(L^{p(\cdot)})} $ and $ g\slash\|g\|_{\ell^{q(\cdot)}(L^{p(\cdot)})} $, as in the proof of Lemma \ref{lem:StrictInequality} there exists $ \{\zeta_\nu\}_\nu $ such that $ \sum_{\nu=1}^{\infty}\zeta_\nu\leq 1 $ and for every $ \nu $ we have
\begin{align}\label{eq:StriclyConvex3}
\varrho_{p(\cdot)}(\frac{f_\nu}{\zeta_\nu^{\frac{1}{q(\cdot)}} 2\mu_1}+\frac{g_\nu}{\zeta_\nu^{\frac{1}{q(\cdot)}} 2\mu_2})<1.
\end{align}
It suffices to show that for every $ \nu $ we have
\begin{align*}
\varrho_{p(\cdot)}(\frac{f_\nu+g_\nu}{\zeta_\nu^{\frac{1}{q(\cdot)}}(\mu_1+\mu_2)})\leq 1.
\end{align*}
Let $ \mu_1\leq\mu_2 $, by reductio and absurdum argument and taking into account Lemma \ref{lem:StrictInequality} we have
\begin{align}\label{eq:Modular and Lemma}
\exists r_1\in(0,\mu_2-\mu_1]~:~1\lneqq\varrho_{p(\cdot)}(\frac{f_\nu}{\zeta_\nu^{\frac{1}{q(\cdot)}}(\mu_1+\mu_2-r_1)}+\frac{g_\nu}{\zeta_\nu^{\frac{1}{q(\cdot)}}(\mu_1+\mu_2+r_1)}). 
\end{align}
Now on the interval $ (0,\mu_2-\mu_1-r_1] $, by using the same procedure as in the Lemma \ref{lem:StrictInequality} for the modular in \eqref{eq:Modular and Lemma}  instead of $ \varrho_{p(\cdot)}(\frac{f_\nu+g_\nu}{\zeta_\nu^{\frac{1}{q(\cdot)}}(\mu_1+\mu_2)}) $, we get into
\begin{align*}
\exists r_2\in(0,\mu_2-\mu_1-r_1]:1\lneqq\varrho_{p(\cdot)}(\frac{f_\nu}{\zeta_\nu^{\frac{1}{q(\cdot)}}(\mu_1+\mu_2-r_1-r_2)}+\frac{g_\nu}{\zeta_\nu^{\frac{1}{q(\cdot)}}(\mu_1+\mu_2+r_1+r_2)}).
\end{align*}
Therefore, by repeating this procedure, there exists a sequence $ \{r_n\}_n $ such that for every $ n\in\mathbb{N} ~:~ r_n\in(0,\mu_2-\mu_1-\sum_{i=1}^{n-1}r_i] $ and
\begin{align*}
1\lneqq\varrho_{p(\cdot)}(\frac{f_\nu}{\zeta_\nu^{\frac{1}{q(\cdot)}}(\mu_1+\mu_2-\sum_{i=1}^{n}r_i)}+\frac{g_\nu}{\zeta_\nu^{\frac{1}{q(\cdot)}}(\mu_1+\mu_2+\sum_{i=1}^{n}r_i)}).
\end{align*}
Passing $ n\rightarrow\infty $, we have $ r_n\rightarrow 0 $ and $ \sum_{i=1}^{n-1}r_i\rightarrow\mu_2-\mu_1 $, thus
\begin{align*}
1\leq\varrho_{p(\cdot)}(\frac{f_\nu}{\zeta_\nu^{\frac{1}{q(\cdot)}}2\mu_1}+\frac{g_\nu}{\zeta_\nu^{\frac{1}{q(\cdot)}}2\mu_2}),
\end{align*}
which in the view of \eqref{eq:StriclyConvex3} is a contradiction and this complete the proof.
\end{proof}


\end{document}